
\input amstex
\documentstyle{amsppt}

\loadeufm  
\loadmsbm  

\def\E{\hbox{${\Cal E}$}}
\def\F{\hbox{${\Cal F}$}}
\def\G{\hbox{${\Cal G}$}}

\def\A{\hbox{${\Cal A}$}}
\def\B{\hbox{${\Cal B}$}}

\def\M{\hbox{${\Cal M}$}}

\def\L{\hbox{${\Cal L}$}}

\def\O#1{\hbox{$\Cal O_{#1}$}}  
\def\co{\hbox{$\Cal O$}}  

\define\SL{\operatorname{SL}}

\magnification \magstep1
\overfullrule=0pt


\topmatter
\title 
Elementary counterexamples to Kodaira \\ vanishing in prime
characteristic
\endtitle

\rightheadtext{Elementary counterexamples to Kodaira}

\author
N.~Lauritzen and A.~P.~Rao
\endauthor
\affil
Aarhus Universitet, \AA rhus, Denmark and 
University of Missouri, St.Louis, USA
\endaffil
\address
Matematisk Institut, Aarhus Universitet, DK-8000 \AA rhus C, Denmark
\endaddress
\email niels\@mi.aau.dk \endemail
\address
University of Missouri -- St.Louis,
St.Louis, MO 63121, USA
\endaddress
\email rao\@arch.umsl.edu \endemail
\keywords Cohomology, line bundles, Kodaira vanishing \endkeywords
\subjclass 14F17, 14M17 \endsubjclass
\abstract
Using methods from the modular representation theory
of algebraic groups one can construct \cite{1} a projective homogeneous 
space for $\SL_4$ in prime characteristic, which violates Kodaira
vanishing. In this note we show how elementary algebraic
geometry can be used to simplify and generalize this example.
\endabstract

\endtopmatter

\document

Let $X$ be a smooth projective variety of dimension $m$ over an 
algebraically closed field of characteristic zero and let \L\ be an ample 
line bundle on $X$. The Kodaira 
Vanishing Theorem \cite{3} states that $H^i(X, \L^{-1}) = 0$ for $i \neq m$. It is well 
known that the result is false in characteristic $p>0$; Raynaud constructed \cite{8}
a smooth projective surface in positive characteristic with an ample line bundle for 
which Kodaira vanishing failed (Mumford \cite{7} had earlier constructed a 
normal non-smooth projective 
surface counter-example). 
Raynaud posed two questions: 
\roster
\item  Are there counter-examples where the vanishing fails for a very
ample line bundle?
\item  Are there pairs $(X,\L)$ ($X$ smooth projective, \L\ ample) for 
which
$$\chi_i(\L\otimes \omega_X) :=h^i(X,\L\otimes \omega_X) - 
h^{i+1}(X,\L\otimes \omega_X) + \dots$$ is not always $\geq 0$?
\endroster

The first author, studying proper homogeneous spaces in characteristic $p$ \cite{4}\cite{5},
answered both questions affirmatively using methods
from modular representation theory of algebraic groups.
The object of this note
is to generalize, via elementary algebraic geometry, the simplest 
counterexample (\cite{1}, Example 4) answering question (1) without using Jantzen's 
sum formula from modular representation theory.
The simplest of our examples
is as follows: let $Y$ be the incidence correspondence of points lying on
planes in projective three space $\Bbb P(V)$.
There is a natural
bundle \G\ of rank $2$ on $Y$ such that the projectivization $\Bbb P(\G)$ 
of $\G$ is the variety of flags in $V$. Now let $X=\Bbb P(F^*\G)$ be the 
projectivization of $F^*\G$ (the Frobenius pull-back of \G.) Then $X$  
can be embedded in $\Bbb P(V) \times \Bbb P(\wedge^2V^\vee) \times 
\Bbb P(V^\vee)$, and the line bundle $\L=\co(1,3,1)$ is very ample on $X$ 
and violates Kodaira vanishing (with $H^5(X,\L^{-1}) \neq 0$). 
\par
It would be interesting to see if such elementary calculations also yield 
examples to Raynaud's second question. The example known so far \cite{5} uses
computer intensive calculations. Also the examples of this note are of
dimension six or more (and with Picard group of rank three), whereas the 
example of Raynaud (and Mumford's normal
variety, earlier) is a surface. While the examples here are more
elementary than Raynaud's example, it would be nice to find other examples
of smaller dimension (or with Picard group of rank two). 
Since it is only the penultimate $H^i(X,\L^{-1})$ which is
non-zero, this failure of Kodaira vanishing is not necessarily inherited
by hyperplane sections.
\par
Both authors would like to thank the Tata Institute of Fundamental Research
for it's hospitality during the period in which this work was done.

\heading 1. Preliminaries \endheading
\subheading{1.1}
Let $V$ be an $n+1$-dimensional vector space over a field $k$.
We will consider $\Bbb P(V)$, $\Bbb P(\wedge^2 V^\vee)$ and $\Bbb P
(V^\vee)$. The tautological line (quotient) bundle of each will be denoted
$\co(1,0,0), \co(0,1,0)$ and $\co(0,0,1)$ respectively. So on $\Bbb P(V)$, there
is an exact sequence of vector bundles
$$ 0 @>>> \A @>>> V \otimes_k \O{\Bbb P(V)} @>>> \co(1,0,0) @>>> 0.$$
This identifies $H^0(\Bbb P(V), \co(1,0,0)$ with $V$.
If we fix a basis $X_0,X_1,\dots ,X_n$ for $V$ and the dual basis $Y_0,Y_1,
\dots , Y_n$ for $V^\vee$, the dual of the above sequence gives the homomorphism
$$ \O{\Bbb P(V)} @>>> \co(1,0,0)\otimes _k V^\vee$$
which determines the global section $\sum X_iY_i$.
\medskip\noindent
\subheading{1.2}
Let $Y$ be a scheme, and
let \E\ be a vector bundle of rank $r$ on $Y$. Let $X=\Bbb P(\E)$ be the projectivized
bundle. It comes with a morphism $\pi: X @>>> Y$, such that $X$ is smooth
over $Y$ and $X$ has a 
tautological line bundle $\O{\pi}(1)$ which appears in a sequence
$$ 0 @>>> \F @>>> \pi^*\E @>>> \O{\pi}(1) @>>> 0.$$
$\omega_{X/Y}$ can be identified with $\wedge^r(\pi^*\E)
\otimes \O{\pi}(-r)$ (\cite{2}, III, Ex 8.4). 
\medskip\noindent
\subheading{1.3}
Let $Y$ be a scheme over a field $k$ of characteristic $p \neq 0$. The 
absolute Frobenius morphism $F:Y@>>> Y$ is defined on the level of affine rings  by mapping
the function $a$ to $a^p$, and has the property that if \L\ is a line bundle
on $Y$, then $F^*\L \cong \L^p$, and if there is a homomorphism $\O X @>>>
\L$ defining a section $s$, it pulls back to a homomorphism  $\O X @>>>
\L^p$ defining the section $s^p$ \cite{6}.
\par
\medskip\noindent
\heading 2. The Example \endheading
Let $V$ be a vector space of dimension $n+1$ over a field $k$ of characteristic
$p \neq 0$ where $n \geq 3$ and $p\geq n-1$. On $\Bbb P(V)$ there is the sequence of bundles
$$  0 @>>> \A @>>> V \otimes_k \O{\Bbb P(V)} @>>> \co(1,0,0) @>>> 0.$$
Let $Y = \Bbb P(\A^\vee)$, with the morphism $\alpha : Y @>>> \Bbb P(V)$.
We have 
$$ 0 @>>> \G @>>> \alpha^*\A^\vee @>>> \O{\alpha}(1) @>>> 0,$$
with \G\ a vector bundle of rank $n-1$. 
The surjection $V^\vee \otimes_k \O{\Bbb P(V)}@>>> \A^\vee$ induces an
inclusion of $\Bbb P(V)$-schemes $ Y \hookrightarrow \Bbb P(V)\times
\Bbb P(V^\vee)$. \par
Let $\pi_1, \pi_2$ be the two projections defined on $\Bbb P(V)\times 
\Bbb P(V^\vee)$. There is the natural map on the product
 $$ \pi_1^*(V^\vee \otimes_k \O{\Bbb P(V)}) @>>> \pi_2^*(\co(0,0,1)) @>>> 0$$
and the composite of this map with the inclusion $ \pi_1^*(\co(-1,0,0)) 
\hookrightarrow \pi_1^*(V^\vee \otimes_k \O{\Bbb P(V)})$ defines a 
homomorphism $\pi_1^*(\co(-1,0,0)) @>>>
\pi_2^*(\co(0,0,1))$.  $Y$ is the zero scheme on the product
$\Bbb P(V)\times \Bbb P(V^\vee)$ of the  global section
of $\co(1,0,1)$ induced by this homomorphism  and
this section can be seen to be just $\sum X_iY_i$. So $Y$ has bihomogeneous coordinate ring
$$k[X_0,\dots,X_n;Y_o,\dots,Y_n]/(\sum X_iY_i)$$ and has canonical line bundle
$\omega_Y = \co(-n,0,-n)$.
\par
If $\beta : Y @>>> \Bbb P(V^\vee)$ is induced by projection $\pi_2$
on the second factor, we see that $ \O{\alpha}(1)$ is the pull-back of
$\co(0,0,1)$. On $Y$ we have the commuting  diagram of exact sequences:

$$ \CD
  @.      0    @.      0  \\
  @.      @AAA       @AAA \\
  0 @>>> \G @>>> \alpha^*\A^\vee @>>> \co(0,0,1) @>>> 0 \\
  @.     @AAA       @AAA                @|             @. \\
  0 @>>> \B @>>>  V^\vee \otimes \O Y @>>> \co(0,0,1) @>>> 0 \\
  @.     @AAA      @AAA                 @.                @. \\
   @.    \co(-1,0,0) @=    \co(-1,0,0) \\
   @.     @AAA                @AAA \\
   @.     0     @.               0 \\
\endCD $$
The middle horizontal sequence is part of the Koszul complex, hence there is
a surjection
$ \wedge^2V^\vee \otimes \co(0,0,-1) @>>> \B @>>> 0$.
\medskip
Now consider $F^*\G'$ where 
$\G' = \G \otimes \co(0,0,1)$. 
Let 
	$X = \Bbb P(F^*\G')$ and let $\pi : X @>>> Y$ be the projection.
$X$ is smooth over $k$ and since there is a surjection $ \wedge^2V^\vee \otimes \O Y
@>>> F^*\G' @>>> 0$, there is an inclusion of $Y$-schemes
   $$ X \hookrightarrow Y \times \Bbb P(\wedge^2V^\vee)$$
It is evident that via projection onto the second factor, $\O{\pi}(1)$
can be identified with the pull-back of $\co(0,1,0)$. Hence the line bundle
$\co(1,1,1)$ on $X$ is very ample. Since $\co(0,n-1,0)$ is globally generated
on $X$,
the line bundle $\L = \co(1,n,1)$ is also very ample on $X$ (\cite{2} II, Ex. 7.5.)
\medskip
\noindent
\proclaim{Claim}: The very ample line bundle \L\ on the smooth variety $X$ 
(smooth over $k$ of dimension $3n-3$)
violates Kodaira vanishing.
\endproclaim
\demo{Proof}: We will compute $H^i(X, \L^{-1})$, where $\L = \co(1,0,1)\otimes
\O{\pi}(n)$. Since $\omega_X = \omega_{X/k} = \omega_{X/Y}\otimes \pi^*(\omega
_Y) = \wedge^{n-1}{F^*\G'}\otimes \O {\pi}(-n+1) \otimes \co(-n,0,-n)
= \co(p,0,p(n-2))\otimes \O {\pi}(-n+1) \otimes \co(-n,0,-n)
= \co (p-n,0,p(n-2)-n) \otimes \O {\pi}(-n+1)$ we get
$$ \align
   H^i(X, \L^{-1}) &= {H^{3n-3-i}(X, \L \otimes \omega_X)}^\vee\\
                &= {H^{3n-3-i}(X, \co(p-n+1,0,p(n-2)-n+1) \otimes \O {\pi}(1))}^\vee\\
                &= {H^{3n-3 -i}(Y, \co(p-n+1, 0,p(n-2)-n+1) \otimes {F^*\G'})}^\vee \\
		&= {H^{3n-3 -i}(Y, \co(p-n+1, 0, p(n-2)-n+1 +p) \otimes {F^*\G})}^\vee\\
		&= {H^{3n-3 -i}(Y, \co(p-n+1, 0, (p-1)(n-1)) \otimes {F^*\G})}^\vee \\
\endalign $$
hence clearly 0 when $3n-3 -i > 2n-1$, the dimension of $Y$, ie. when $i <
n-2$.\par
Let $\M = \co(p-n+1, 0, (p-1)(n-1))$ on $Y$. The Frobenius pull-back of the 
commuting diagram above, tensored with \M\ gives

$$ \CD
  @.      0    @.      @.  \\
  @.      @AAA       @. \\
 @.       \M\otimes F^*\G  @.  @. \\
 @.       @AAA       @. \\
 0 @>>>  \M\otimes F^*\B @>>> V^\vee \otimes \M @>>> \M\otimes \co(0,0,p)@>>> 0.\\
@.       @AAA        @. \\
@.       \M\otimes \co(-p,0,0) @. \\
@.      @AAA       @. \\
@.       0       @. \\
\endCD $$

By running along the exact sequence on $\Bbb P(V)\times \Bbb P(V^\vee)$
$$ 0 @>>> \co(-1,0,-1) @>>> \co @>>> \O Y @>>> 0, $$
we see that $H^j(Y, \M\otimes \co(-p,0,0)) = H^j(Y, \co(1-n,0,(p-1)(n-1)))$ 
is always zero, hence
$$  H^i(X, \L^{-1}) = {H^{3n-3 -i}(Y, \M \otimes F^*\B )}^\vee. $$
$V^\vee \otimes \M$ and $\M\otimes \co(0,0,p)$ have only $H^0$ as nonzero
cohomology ($H^0$ is non-zero since $p \geq n-1$)
hence $F^*\B \otimes \co(p+1-n, 0, (p-1)(n-1))$ has zero cohomology except
possibly $H^0$ and
$H^1$ which are to be studied by considering the map
$$ A: V^\vee \otimes H^0(Y, \co(p+1-n, 0, (p-1)(n-1))) @>>> 
   H^0(Y, \co(p+1-n, 0, (p-1)(n-1) + p)).$$
The homomorphism $A$ is given by the matrix $[Y_0^p,Y_1^p,\dots,Y_n^p]$ and it
is easy to see that it is not onto: Consider the element
$$ t = X_n^{p+1-n}.Y_0Y_1^{p-1}Y_2^{p-1}\dots Y_n^{p-1}$$
which is well defined in $H^0(Y, \co(p+1-n, 0, (p-1)(n-1) + p))$ modulo multiples
of $X_0Y_0+X_1Y_1+\dots + X_nY_n$. No element in its equivalence class can be
expressed as a sum of multiples of $Y_0^p,Y_1^p,\dots,Y_n^p$, ie. $t$ is not
in the image of $A$.
\par
It follows that with $3n-3-i = 1$, $H^i(X, \L^{-1})\neq 0$. Hence 
$H^{3n-4}(X, \L^{-1}) \neq 0$, where $X$ has dimension $3n-3$  and $H^i(
X, \L^{-1}) = 0$ for all $i < 3n-4$.\qed 
\enddemo
\vfill\eject

\enddocument